\documentclass[12pt]{amsart}

\usepackage{verbatim}
\usepackage{amssymb}
\usepackage{upref}
\usepackage[all]{xy}

\newtheorem{thm}{Theorem}[section]
\newtheorem{lem}[thm]{Lemma}
\newtheorem{cor}[thm]{Corollary}
\newtheorem{prop}[thm]{Proposition}

\theoremstyle{definition}
\newtheorem{defn}[thm]{Definition}

\newtheorem*{notn}{Notation and Terminology}

\newtheorem{rem}[thm]{Remark}

\numberwithin{equation}{section}

\newcommand{\secref}[1]{Section~\textup{\ref{#1}}}

\newcommand{\thmref}[1]{Theorem~\textup{\ref{#1}}}
\newcommand{\corref}[1]{Corollary~\textup{\ref{#1}}}
\newcommand{\lemref}[1]{Lemma~\textup{\ref{#1}}}
\newcommand{\propref}[1]{Proposition~\textup{\ref{#1}}}

\newcommand{\remref}[1]{Remark~\textup{\ref{#1}}}

\newcommand{\diagref}[1]{Diagram~\textup{\ref{#1}}}

\renewcommand{\labelenumi}{(\roman{enumi})}

\newcommand{\midtext}[1]{\quad\text{#1}\quad}
\newcommand{\righttext}[1]{\qquad\text{#1 }}
\renewcommand{\and}{\midtext{and}}
\renewcommand{\for}{\righttext{for}}
\newcommand{\all}{\righttext{for all}}

\newcommand{\textin}{\midtext{in}}

\newcommand{\C}{\mathbb C}

\newcommand{\CC}{\mathcal C}
\newcommand{\KK}{\mathcal K}

\renewcommand{\AA}{\mathcal A}

\renewcommand{\a}{\alpha}
\renewcommand{\b}{\beta}
\renewcommand{\d}{\delta}
\renewcommand{\l}{\lambda}

\newcommand{\p}{\phi}
\newcommand{\e}{\epsilon}
\newcommand{\x}{\xi}
\newcommand{\m}{\mu}

\newcommand{\s}{\sigma}
\newcommand{\g}{\gamma}
\renewcommand{\t}{\theta}

\renewcommand{\epsilon}{\varepsilon}

\renewcommand{\P}{\Phi}

\renewcommand{\L}{\Lambda}

\DeclareMathOperator{\aut}{Aut}
\DeclareMathOperator{\ad}{Ad}

\DeclareMathOperator{\obj}{Obj}

\DeclareMathOperator*{\spn}{span}
\DeclareMathOperator*{\clspn}{\overline{\spn}}

\newcommand{\id}{\text{\textup{id}}}

\newcommand{\<}{\langle}
\renewcommand{\>}{\rangle}

\newcommand{\inv}{^{-1}}

\newcommand{\iso}{\overset{\cong}{\longrightarrow}}
\renewcommand{\bar}{\overline}
\newcommand{\what}{\widehat}
\newcommand{\wilde}{\widetilde}

\newcommand{\rt}{\textup{rt}}

\newcommand{\dn}{\downarrow}

\begin{document}

\title{Categorical Landstad duality for actions}

\author{S. Kaliszewski}

\address{Department of Mathematics and Statistics
\\Arizona State University
\\Tempe, Arizona 85287}

\email{kaliszewski@asu.edu}

\author{John Quigg}

\address{Department of Mathematics and Statistics
\\Arizona State University
\\Tempe, Arizona 85287}

\email{quigg@asu.edu}



\date{November 12, 2007}

\begin{abstract}
We show that the category $\AA(G)$ of actions of
a locally compact group
$G$ on $C^*$-algebras
(with equivariant nondegenerate $*$-homomorphisms into multiplier
algebras) is equivalent, via a full-crossed-product functor,
to a comma category of maximal coactions of $G$
under the comultiplication $(C^*(G),\delta_G)$; and
also that $\AA(G)$ is equivalent, via a reduced-crossed-product
functor,
to a comma category of normal coactions under the comultiplication.
This extends classical Landstad duality to a category equivalence,
and allows us to identify those $C^*$-algebras which are isomorphic to
crossed products by $G$ as precisely those which form part of an
object
in the appropriate comma category.
\end{abstract}

\maketitle

\section*{Introduction}
\label{intro}
\emph{Landstad duality} (a term coined by the second author in
\cite{Q92}) refers to a particular method of characterizing
crossed-product $C^*$-algebras.
The first appearance of this method is in
\cite{lan:dual}, where Landstad
characterized reduced crossed products by actions of locally compact
groups in terms of the existence of suitably equivariant reduced
coactions.
The second author proved a dual version
of this
in \cite{Q92},
giving a characterization of crossed products by coactions in terms
of the existence of suitably equivariant actions.
In \cite{KQ07}, the authors applied the recently-developed
theory of maximal coactions (see \cite{ekq}) to give a version of
Landstad's characterization for \emph{full}, rather than reduced,
crossed products by actions.

In the present paper we analyze the method of Landstad duality more
deeply, shifting the focus from characterizing crossed products to
developing a process which recovers the action
\emph{up to isomorphism}
from its crossed
product.
It is useful to compare this with the more
well-established \emph{crossed-product duality}, which uses the
crossed product by the
dual coaction to recover the action up to Morita-Rieffel equivalence
(see \cite[Appendix]{enchilada}
for a recent survey).

One of our objectives in the present paper is to promulgate a
``categorical imperative'': in order to have a robust and complete
theory, it is not enough to establish results for the $C^*$-algebras
alone --- rather, one must also take care of the morphisms. There are
numerous results in the literature concerning equivariant homomorphisms
(and isomorphisms in particular), and categorical techniques give a
unified way of dealing with them. Together with Echterhoff and
Raeburn, we established such a framework in \cite{enchilada} for
$C^*$-correspondences, and in the present paper we
do it for homomorphisms. One of the main differences between the
categories of \cite{enchilada} and the present paper is what ``isomorphism'' means --- Morita-Rieffel equivalence in the memoir and ordinary isomorphism here.

Our two versions of categorical Landstad duality
(Theorem~\ref{actionnormal} for reduced crossed products and Theorem~\ref{full landstad} for full crossed products) give equivalences
between the category of actions of a given group $G$ and the
categories of normal and maximal coactions, respectively,
equipped with suitably equivariant
homomorphisms of $G$
(see \secref{prelim} for the definitions).
The power of the categorical approach is manifested via the use of
\emph{comma categories} to encode the ``suitably equivariant
homomorphisms''.

We give a third category equivalence
(Theorems~\ref{maximalnormalnocomma} and \ref{maximalnormal} --- the
second includes appropriate comma categories), between maximal and
normal coactions. In fact, we use this to deduce our categorical
Landstad duality for full crossed products from the version for
reduced crossed products.
The method of proof of \thmref{maximalnormalnocomma} deserves
comment: maximal and normal coactions of $G$ form subcategories of
the category of all $G$-coactions, and they are related to the
ambient category in ``dual'' ways. Every coaction is a quotient of a
maximal one, and has a normal one as a quotient, and in both cases
the quotient homomorphisms give isomorphic crossed products.
Our results in \secref{maximal normal} are applications of a study of
this situation in an abstract setting:
Propositions~\ref{subcategory} and \ref{comma}
are two general results giving category equivalences
which may be well-known to category \emph{cognoscenti},
although we could not find them in the literature.

The results on Landstad duality currently in the literature lead one to
\emph{suspect} that recovering an action from the crossed product and a
suitable equivariant coaction is some sort of ``inverse process'', but
it is only through the use of categorical methods that this
suspicion is fully justified and made precise.

Similarly, the universal properties of maximal and normal coactions
lead one to \emph{suspect} that these two types of coaction are
``essentially the same'', but again it is the categorical framework of
the present paper that tells the full story.

We hope it will become clear that
the results of the present paper
do not consist of merely translating pre-existing results on
Landstad duality into categorical language.
For example, handling the morphisms in
the proof of
\thmref{actionnormal}
requires new arguments, which
occupy four pages.

The results of this paper are intended primarily for the use of operator algebraists.
The theory of noncommutative duality for $C^*$-dynamical systems plays a central role, but we have tried to make our techniques accessible to operator algebraists who might not be experts in noncommutative duality.
Our main results involve category theory in crucial ways,
so we emphasize categorical techniques throughout.
For this reason, we give a somewhat more detailed description of the preliminary material than is typical in a paper on noncommutative duality.

The  authors thank Iain Raeburn for helpful conversations.

\section{Preliminaries}
\label{prelim}

\subsection*{Category theory}

We adopt the conventions of \cite{maclane} for category theory,
with \cite{awodey} as a secondary reference.
If $C$ is a category then $\obj C$ denotes the class of objects.
If $x,y\in \obj C$ then $C(x,y)$ denotes the set of arrows of $C$ with
domain $x$ and codomain $y$.
We write
``$f:x\to y$ in $C$'' to mean $f\in C(x,y)$.

Recall the following (more-or-less) standard definitions: if $F:C\to D$ is a functor, then
\begin{enumerate}
\item $F$ is \emph{full} (respectively, \emph{faithful})
if it maps $C(x,y)$ surjectively (respectively, injectively) to $D(Fx,Fy)$
for all $x,y\in \obj C$;
\item a functor $G:D\to C$
is a \emph{quasi-inverse} of $F$ if
$FG\cong \id_D$ and $GF\cong \id_C$ (where ``$\cong$'' here means ``naturally equivalent'' and ``\id'' refers to the identity functor);
\item $F$ is an \emph{equivalence} if it has a quasi-inverse
(in the sense of (ii));
\item $F$ is \emph{essentially surjective} if every object in $D$ is isomorphic (in $D$) to one in the image of $F$.
\end{enumerate}
Of the above definitions, only (i) and (iii) appear in \cite{maclane}, while all appear in \cite{awodey}.
Warning: Mac Lane \cite[discussion preceding Section~IV.1, Corollary~1]{maclane} uses the term \emph{quasi-inverse} in a somewhat different sense.

Of course, a functor $F:C\to D$ is an \emph{isomorphism} if it is full, faithful, and bijective on objects, and then the inverse map $G:D\to C$ is a functor for which $FG=\id_D$ and $GF=\id_C$.

The following result is standard in category theory; we briefly outline the proof for convenient reference in the proof of the subsequent corollary.

\begin{prop}
[{\cite[Section~IV.4, Theorem 1]{maclane}}]
\label{equivalence}
A functor $F:C\to
D$ is an equivalence if and only if it is full, faithful, and essentially surjective.
\end{prop}

\begin{proof}[Outline of proof]
The nontrivial direction is the reverse implication: assuming $F$ is full, faithful, and essentially surjective,
for each $u\in \obj D$
choose
(using an extension of the Axiom of Choice for classes)
an object $Gu$ of $C$
and an isomorphism $\t_u:FGu\iso u$ in $D$.
For each $f\in D(u,v)$, let $Gf\in C(Gu,Gv)$ be the unique arrow such that
\[
FGf=\t_v\inv\circ f\circ \t_u.
\]
Then $G$ is a quasi-inverse of $F$.
\end{proof}

\begin{cor}
\label{equivalence2}
If
$F:C\to D$ and $G:D\to C$
are functors such that
$FG\cong \id_D$ and
$F$ is full and faithful,
then $G$ is a quasi-inverse of $F$.
\end{cor}

\begin{proof}
Let $\t:FG\iso \id_D$.
Then for $f\in D(u,v)$ we have $FGf=\t_v\inv\circ f\circ \t_u$, so $G$ coincides with the functor constructed in the proof \propref{equivalence},
and in that proof $G$
is shown to be a quasi-inverse of~$F$.
\end{proof}


Recall \cite[Section~IV.1]{maclane}
that an \emph{adjunction} $(F,G,\p):C\to D$ comprises functors $F:C\to D$ and $G:D\to C$, and for all $x\in \obj C$ and $y\in \obj D$ a bijection $\p=\p_{x,y}:D(Fx,y)\to C(x,Gy)$ which is natural in $x,y$. $F$ is called a \emph{left adjoint} of $G$ and $G$ is called a \emph{right adjoint} of $F$.
The adjunction is uniquely determined by its \emph{unit}, which is a natural transformation $\eta:\id_C\to GF$ such that for all $x\in \obj C$ the map $\eta_x:x\to GFx$ is \emph{universal from $x$ to $G$}, meaning that for all $y\in \obj D$ and $f\in C(x,Gy)$ there is a unique $g\in D(Fx,y)$ (namely, $g=\p\inv(f)$) such that the triangle in the following diagram commutes in $C$:
\[
\xymatrix{
x \ar[r]^-{\eta_x} \ar[dr]_f
&GFx \ar@{-->}[d]^{Gg}
&Fx \ar@{-->}[d]^g_{!}
\\
&Gy
&y.
}
\]
In fact, the adjunction is completely determined by the functor $G:D\to C$ and, for each $x\in \obj C$, a universal arrow $\eta_x$ from $x$ to $G$.
Dually, the adjunction is uniquely determined by its \emph{counit}, which is a a natural transformation $\e:FG\to \id_D$ such that for all $y\in \obj D$ the map $\e_y:FGy\to y$ is \emph{universal from $F$ to $y$}, meaning that for all $x\in \obj C$ and $g\in D(Fx,y)$ there is a unique $f\in C(x,Gy)$ (namely, $f=\p(g)$) such that the triangle in the following diagram commutes in $D$:
\[
\xymatrix{
Gy
&y
&FGy \ar[l]_-{\e_y}
\\
x \ar@{-->}[u]^f_{!}
&&Fx \ar@{-->}[u]_{Ff} \ar[ul]^g.
}
\]


We will make extensive use of comma categories, and we record our
notation for them here.

\begin{defn}
Let $a$ be an object in a category $C$.
The \emph{comma category of objects of $C$
under $a$}, written $a\downarrow C$, has objects $(x,f)$, where $f:a\to x$ in $C$; an arrow
$h:(x,f)\to (y,g)$ in $a\dn C$ is an arrow $h:x\to y$ in $C$ such that $h\circ
f=g$.
\end{defn}

Heuristically, an object in $a\downarrow C$ is an arrow
\[\xymatrix{a\ar[d]_f\\x,}\]
and an arrow in $a\downarrow C$ is a commuting triangle
\[\xymatrix{&a\ar[dl]_f\ar[dr]^g\\x\ar[rr]_-h&&y.}\]

\begin{defn}
If $D$ is a subcategory of $C$, it is easy to check that we get a subcategory
$a\dn D$
of $a\dn C$
by restricting the $x$, $y$, and $h$ (but
\emph{not} the $f$ or $g$!) to be in $D$.
\end{defn}

\begin{defn}
Dually, the \emph{comma category of objects of $C$ over $a$}, written $C\downarrow a$,
has objects
$(x,f)$ where $f:x\to a$ in $C$ and arrows $h:(x,f)\to (y,g)$ where
$h:x\to y$ in $C$ and $g\circ h=f$:
\[\xymatrix{x\ar[rr]^h\ar[dr]_f&&y\ar[dl]^g\\&a.}\]
Again, if $D$ is a subcategory of $C$ we can form the subcategory
$D\downarrow a$ of $C\downarrow a$.
\end{defn}

\subsection*{Fundamental category}

Our fundamental category, denoted $\CC$, has the class of all $C^*$-algebras as its objects. An arrow $\p:A\to B$ in the category $\CC$ is by definition a nondegenerate homomorphism $\p:A\to M(B)$
--- this works because such a homomorphism extends uniquely to a unital homomorphism $\p:M(A)\to M(B)$, and hence nondegenerate homomorphisms can be composed.
A key feature of this category is that isomorphisms coincide with ordinary isomorphisms of $C^*$-algebras.
Warning: \cite{enchilada} uses the same
notation for a category whose arrows are isomorphism classes of
$C^*$-correspondences (called
right-Hilbert bimodules in \cite{enchilada}), but
in this paper we will have no use for correspondences as arrows.

\subsection*{Actions and coactions}

We adopt, more or less, the conventions of \cite{enchilada},
\cite{ekq},
and \cite{Q94} for actions and coactions
of locally compact Hausdorff groups on $C^*$-algebras.
An \emph{action} of a locally compact group $G$ on a $C^*$-algebra $A$ is a homomorphism $\a$ from $G$ to the automorphism group $\aut A$ such that $s\mapsto \a_s(a)$ is continuous for each $a\in A$.
A \emph{coaction} of $G$ on $A$ is a nondegenerate homomorphism $\d:A\to M(A\otimes C^*(G))$ (where we use the minimal $C^*$-tensor product throughout) such that $\clspn\{\d(A)(1\otimes C^*(G))\}=A\otimes C^*(G)$ and $(\d\otimes\id)\circ\d=(\id\otimes\d_G)\circ\d$,
where $\d_G:C^*(G)\to M(C^*(G)\otimes C^*(G))$ is the unique homomorphism such that $\d_G(s)=s\otimes s$ for all $s\in G$ (after extending $\d_G$ to the multiplier algebra).
$\d_G$ is called the \emph{canonical coaction} of $G$ on $C^*(G)$, or the \emph{comultiplication} on $C^*(G)$.

$\AA(G)$
denotes the category
in which an object is an action $(A,\a)$ of $G$,
and in which an arrow $\p:(A,\a)\to (B,\b)$ is an arrow $\p:A\to B$ in
$\CC$ which is $\a-\b$ equivariant.
Similarly, $\CC(G)$ denotes the category
in which an object is a coaction $(A,\d)$, and in which an arrow $\p:(A,\d)\to (B,\e)$ is an arrow
$\p:A\to B$ in $\CC$ which is $\d-\e$ equivariant.

\begin{notn}
For a coaction $(A,\d)$,
although an object in the comma category $(A,\d)\downarrow
\CC(G)$ is officially an ordered pair $((B,\e),\p)$ with
$\p:(A,\d)\to (B,\e)$ in $\CC(G)$, we find it convenient to
write this as a triple $(B,\e,\p)$, and we also let
$(A,\d,\p)$
denote the
corresponding object in
the comma
category $\CC(G)\downarrow (B,\e)$.
\end{notn}

We freely identify an arrow $u:C^*(G)\to A$ in $\CC$
with a strictly continuous unitary homomorphism $u:G\to M(A)$.
Every such $u$ implements an \emph{inner action} $\ad u$ of $G$ on $A$,
and this determines the object map
\[
(A,u)\mapsto (A,\ad u)
\]
of a functor
\[
\ad:C^*(G)\dn \CC\to \AA(G).
\]
An arrow $\p:(A,u)\to (B,v)$ in the comma category $C^*(G)\dn \CC$ is an arrow $\p:A\to B$ in $\CC$ such that $\p\circ u=v$
(extending $\p$ to $M(A)$ in the usual manner),
and the functor $\ad$ takes this to the arrow $\p:(A,\ad u)\to (B,\ad v)$.

The crossed-product construction for actions can be regarded as a left adjoint of the functor $\ad$.
Operator algebraists are accustomed to
giving a particular
construction of the crossed product, then verifying that
it
is functorial from $\AA(G)$ to $\CC$ (that is, only keeping the
$C^*$-algebra from the covariant homomorphism)
by applying the universal property.
The viewpoint we advocate in this paper is that there is no need for us to make a
choice of the crossed product --- universality makes the
crossed product functorial no matter what choice is made, and any
other
choice of crossed product would give a naturally isomorphic functor.
Because this idea is central for the various applications we make
in this paper, for clarity we give more detail in this case.
The idea is to choose crossed products arbitrarily
(appealing to the Axiom of Choice for classes),
and then
universality tells us what to do for the
arrows.
More precisely,
for each action $(A,\a)$ there is a universal arrow
from $(A,\a)$ to the functor $\ad$.
This means (see \cite[Section~IV.1]{maclane}) that there is an object,
which in the standard notation of operator algebras would be denoted by
$(A\times_\a G,i_G)$,
of the comma category $C^*(G)\dn \AA(G)$,
and an arrow,
which operator algebraists would denote by $i_A$,
from $(A,\a)$ to the associated inner action $(A\times_\a G,\ad i_G)$,
such that for every object $(B,u)$ of $C^*(G)\dn \CC$ and every arrow $\pi:(A,\a)\to (B,\ad u)$ in $\AA(G)$ there is a unique arrow,
which operator algebraists would denote by $\pi\times u:(A\times_\a G,i_G)\to (B,u)$,
such that the triangle in the
following diagram commutes in $\AA(G)$:
\[
\xymatrix@C+30pt{
(A,\a) \ar[r]^-{i_A} \ar[dr]_\pi
&(A\times_\a G,\ad i_G) \ar@{-->}[d]^{\pi\times u}
&(A\times_\a G,i_G) \ar@{-->}[d]^{\pi\times u}_{!}
\\
&(B,\ad u)
&(B,u).
}
\]
By \cite[Section~IV.1, Theorem~2]{maclane}, the existence for each action $(A,\a)$ of a universal arrow $i_A$ guarantees that the assignment $(A,\a)\mapsto (A\times_\a G,i_G)$ is the object map of a functor $\times G:\AA(G)\to C^*(G)\dn \CC$ which is a left adjoint of the functor $\ad$, and moreover (by \cite[Section~IV.1, Corollary~2]{maclane}) that any other choice of universal arrows would give a naturally isomorphic functor.
The functor $\times G$ sends an arrow $\p:(A,\a)\to (B,\b)$ in $\AA(G)$ to an arrow, denoted by operator algebraists as $\p\times G$, from $(A\times_\a G,i_G)$ to $(B\times_\b G,i_G)$.
The \emph{crossed product} of the action $(A,\a)$ is usually denoted by the triple $(A\times_\a G,i_A,i_G)$, and is characterized by the above universal property.
If $\pi:(A,\a)\to (B,\ad u)$ in $\AA(G)$,
the pair $(\pi,u)$ is usually called a \emph{covariant homomorphism} of the action $(A,\a)$, and the associated arrow $\pi\times u:A\times_\a G\to B$ is called the \emph{integrated form} of $(\pi,u)$.
There is an obvious category whose objects are triples $(B,\pi,u)$ (with $B,\pi,u$ as above), and in which an arrow $\s:(B,\pi,u)\to (C,\p,v)$ consists of an arrow $\s:B\to C$ in $\CC$ which is simultaneously an arrow
$(B,\pi)\to (C,\p)$ in $A\dn \CC$ and an arrow
$(B,u)\to (C,v)$ in $C^*(G)\dn \CC$.
The universal property of crossed products says that
$(A\times_\a G,i_A,i_G)$ is an initial object in this category. We will find it convenient to use the crossed-product functor in the form $(A,\a)\mapsto (A\times_\a G,i_A,i_G)$.

Thus, we can, and will, assume that some choice of crossed-product
functor has been made, but it is completely irrelevant for our
purposes
which particular choice this is.

A completely parallel treatment exists for for coaction-crossed products:
an object $(A,\m)$ of $C_0(G)\dn \CC$
implements an \emph{inner coaction} $\ad \m$ of $G$ on $A$,
and this determines the object map $(A,\m)\mapsto (A,\ad \m)$ of a functor $\ad:C_0(G)\dn \CC\to \CC(G)$.
For each coaction $(A,\d)$ there are
an object $(A\times_\d G,j_G)$ of $C_0(G)\dn \CC$ and
a universal arrow
\[
j_A:(A,\d)\to (A\times_\d G,\ad j_G)
\]
from $(A,\d)$ to the functor $\ad$.
Universality means that for every object $(B,\m)$ of $C_0(G)\dn\CC$ and every arrow $\pi:(A,\d)\to (B,\ad\m)$ in $\CC(G)$ there is a unique arrow $\pi\times\m:(A\times_\d G,j_G)\to (B,\m)$ in $C_0(G)\dn\CC$ such that the triangle in the following diagram commutes in $\CC(G)$:
\[
\xymatrix@C+30pt{
(A,\d) \ar[r]^-{j_A} \ar[dr]_\pi
&(A\times_\d G,\ad j_G) \ar@{-->}[d]^{\pi\times \m}
&(A\times_\d G,j_G) \ar@{-->}[d]^{\pi\times \m}_{!}
\\
&(B,\ad \m)
&(B,\m).
}
\]
Universality gives a left adjoint functor
$\times G:\CC(G)\to C_0(G)\dn \CC$ of $\ad$,
with object map $(A,\d)\mapsto (A\times_\d G,j_G)$.
The image of an arrow $\p:(A,\d)\to (B,\e)$ in $\CC(G)$ under the functor $\times G$ is denoted $\p\times G$.
The crossed product of the coaction $(A,\d)$ is usually denoted
$(A\times_\d G,j_A,j_G)$, and is an initial object of an appropriate category of covariant homomorphisms.

For both the above action- and coaction-product functors, the various parts of the associated adjunctions can be identified in terms of standard crossed-product concepts. For example, in the case of actions, the unit is
\[
\eta_{(A,\a)}=i_A:(A,\a)\to (A\times_\a G,\ad i_G)
\]
and the co-unit is 
\[
\e_{(A,u)}=\id_A\times u:(A\times_{\ad u} G,i_G).
\]
Since we do not need to work with these assorted bits of the adjunctions in this paper, we refrain from going into more detail.

If $(A,\a)$ is an action, then $\what\a$ denotes
the dual coaction of $G$ on $A\times_\a G$,
determined by the covariant homomorphism
\[
(i_A\otimes 1,(i_G\otimes \id)\circ \d_G).
\]
A fundamental property of the dual coaction is that
the arrow $i_G:C^*(G)\to A\times_\a G$ is also
an arrow in $\CC(G)$:
\[i_G:(C^*(G),\d_G)\to (A\times_\a G,\what\a).\]
We can, and will, also regard the crossed product as a functor
\[
(A,\a)\mapsto (A\times_\a G,\what\a,i_G)
\]
from $\AA(G)$ to $(C^*(G),\d_G)\dn \CC(G)$.

Similarly, if
$(A,\d)$ is a coaction, then $\what\d$ is
the dual action of $G$ on $A\times_\d G$,
where for each $s\in G$ the automorphism $\what\d_s$ is
determined by the covariant homomorphism
\[
(j_A,j_G\circ \rt_s),
\]
where in turn $\rt$ denotes the action of $G$ on $C_0(G)$ given by right translation:
\[
\rt_s(f)(t)=f(ts)\righttext{for}f\in C_0(G).
\]
We can also regard the crossed product as a functor
\[
(A,\d)\mapsto (A\times_\d G,\what\d)
\]
from $\CC(G)$ to $\AA(G)$.\
\footnote{Although it could be regarded as a functor into an appropriate comma category, as we did for the action-crossed product, we will have no
need of that for coaction-crossed products in this paper.}

\subsection*{Maximal and normal coactions}

A coaction $(A,\d)$ is \emph{maximal} \cite{ekq}
if full-crossed-product duality holds, that is, if
the canonical
surjection
\[\Phi=\Phi_{(A,\d)}:A\times_\d G\times_{\what\d} G
\to A\otimes \KK(L^2(G))\]
is an isomorphism,
where $\KK(L^2(G))$ denotes the $C^*$-algebra of compact linear operators on the Hilbert space $L^2(G)$,
and is \emph{normal} \cite{Q94} if the arrow
\[j_A:A\to A\times_\d G\textin \CC\]
is injective\footnote{and there is also a characterization in terms of
reduced-crossed-product duality \cite[Proposition~2.2]{ekq}, but we
will not need it in this paper}.
Let $\CC^m(G)$ and $\CC^n(G)$ denote the full subcategories of
$\CC(G)$ whose objects are
maximal and normal coactions, respectively.

A \emph{maximalization} of a coaction $(A,\d)$ is a surjection $\psi:(B,\e)\to (A,\d)$ such that $(B,\e)$ is maximal and $\psi\times G:B\times_\e G\to A\times_\d G$ is an isomorphism, and $\psi$ is called a \emph{maximalizing map}.
A \emph{normalization} of $(A,\d)$ is a surjection $\eta:(A,\d)\to (B,\e)$ such that $(B,\e)$ is normal and $\eta\times G:A\times_\d G\to B\times_\e G$ is an isomorphism, and $\eta$ is called a \emph{normalizing map}.
Maximalizations and normalizations always exist
(see \cite[Theorem~3.3]{ekq} for maximalizations and
\cite[Proposition~2.3]{Q94} for normalizations).
The requirement that maximalizing and normalizing maps be surjective is redundant, as \lemref{onto} below shows.

For an action $(A,\a)$
we let $\L=\L_{(A,\a)}:A\times_\a G\to A\times_{\a,r} G$ denote a
choice of
regular representation onto (a choice of) the reduced crossed
product, and in particular, taking $A=\C$, we let
$\l:C^*(G)\to C^*_r(G)$ denote the regular representation onto
the reduced $C^*$-algebra of $G$.
It does not matter how the reduced crossed product is chosen, as long
as the kernel of $\L$ in $A\times_\a G$ is the correct ideal. The
surjection $\L$ is the integrated form of a covariant homomorphism
which
we denote by $(i^r_A,i^r_G)$.
There is a unique coaction\footnote{Warning: the dual coaction on the reduced crossed product has been denoted by $\what\a^n$ (for example) in other papers --- the notation $\what\a^r$ used here is new; we introduce it in part to reduce conflicts with other coaction constructions (for example, normalization and reduction).}
$\what\a^r$
such that
\[\L:(A\times_\a G,\what\a)\to (A\times_{\a,r} G,\what\a^r)
\textin \CC(G),\]
and in particular there is a unique coaction $\d_G^r$ such that
\[\l:(C^*(G),\d_G)\to (C^*_r(G),\d_G^r).\]
Given an arrow $\p:(A,\a)\to (B,\b)$ in $\AA(G)$,
there is a unique arrow $\p\times_r G$ in $\CC(G)$ such that
\[\xymatrix@C+20pt{
(A\times_\a G,\what\a) \ar[r]^-{\L_{(A,\a)}} \ar[d]_{\p\times G}
&(A\times_{\a,r} G,\what\a^r) \ar@{-->}[d]^{\p\times_r G}_{!}
\\
(B\times_\b G,\what\b) \ar[r]_-{\L_{(B,\b)}}
&(B\times_{\b,r} G,\what\b^n)
}\]
commutes.
The associated reduced-crossed-product functor
from $\AA(G)$ to $\CC(G)$
is given by
\[
(A,\a)\mapsto (A\times_{\a,r} G,\what\a^r)
\and
\p\mapsto \p\times_r G,
\]
and we can also regard it as a functor
\[
(A,\a)\mapsto (A\times_{\a,r} G,\what\a^r,i^r_G).
\]
from $\AA(G)$ to $(C^*(G),\d_G)\dn \CC(G)$.
Moreover, by construction the arrows $\L_{(A,\a)}$ give a natural
transformation
between the full- and reduced-crossed-product functors.

The dual coaction $(A\times_\a G,\what\a)$ is maximal
(see \cite[Proposition~3.4]{ekq} and \cite[Theorem~3.7]{ekr}),
and
$(A\times_{\a,r} G,\what\a^r,\L_{(A,\a)})$ is a normalization of
it (\cite[Proposition~A.61]{enchilada}).

\label{boilerplate}
We occasionally need to appeal to results from the literature concerning
\emph{reduced} coactions, and apply them to \emph{normal} coactions.
This is justified by the following facts,
collected from \cite[Proposition~3.3, Proposition~3.8,
and Theorem~4.7]{Q94} (see also \cite{rae:representation}):
\begin{itemize}
\item
if $(A,\d)$
is
a normal coaction, then the homomorphism\footnote{The notation $\d_\l$ we introduce here is new, and is meant to emphasize the dependence upon the choice of the regular representation $\l$. The notation $\d^r$ commonly appears in the literature.}
$\d_\l$
defined by the
commuting diagram
\[\xymatrix{
A\ar[r]^-\d\ar[dr]_{\d_\l}&M(A\otimes C^*(G))\ar[d]^{\id\otimes\l}
\\&M(A\otimes C^*_r(G))
}\]
is a reduced coaction of $G$ on $A$;

\item
$(A\times_{\d_\l} G,\d_\l,1\otimes M)$
is a crossed product of $(A,\d)$;

\item
the dual action on $A\times_{\d_\l} G$ is defined in the same way for $\d$ and for $\d_\l$;

\item
if $(A,\e)$ is a reduced coaction, then there exists a
unique normal coaction $(A,\d)$ such that $\e=\d_\l$.
\end{itemize}

\section{Elementary category theory results}
\label{nonsense}

In this section we record a couple of general results in category theory
which may be well-known, but since we could not find them in the
literature we give the proofs for completeness.
First we prove an abstract category-equivalence result for
subcategories, which we then promote to an equivalence between comma
categories.
We will later apply these results to obtain an equivalence between
maximal and normal coactions.

Throughout this section,
we consider the following situation: $D$ and
$E$ are subcategories of a category $C$,
with $D$ coreflective and $E$
reflective.

Recall
from \cite[Section~IV.3]{maclane}
that $E$ is \emph{reflective} in $C$ if
the inclusion functor $E\hookrightarrow C$ has a left adjoint
$F:C\to E$, called the \emph{reflector}.
The adjunction is called a \emph{reflection} of $C$ in $E$.
By the general theory of adjoints,
$E$ is reflective in $C$ if and only if
for each $x\in \obj C$ the comma category $x\downarrow E$ has an initial
object.
Recall that
we defined $x\dn E$ as a subcategory of $x\dn C$.
Thus,
letting $\eta$ denote the unit of the reflection,
for each $x\in \obj C$ and $(y,f)\in x\dn E$ there is a unique $g\in E(Fx,y)$ making the diagram
\[
\xymatrix{
x \ar[d]_{\eta_x} \ar[dr]^f
\\
Fx \ar@{-->}[r]_g^{!}
&y
}
\]
commute in $C$; we emphasize that $Fx,y,g$ are in the subcategory $E$, but $x,\eta_x,f$ are in the ambient category $C$.

Similarly, \emph{coreflectivity} of $D$ means that
the inclusion functor $D\hookrightarrow C$ has a right adjoint $G:C\to D$,
called the \emph{coreflector}, and the adjunction is called a \emph{coreflection} of $C$ in $D$.
Thus, letting $\psi$ denote the counit of the coreflection,
for all
$x\in \obj C$ and
$(y,f)\in D\dn x$
there exists a unique $g\in D(y,Gx)$
making the diagram
\[\xymatrix{y\ar[dr]_f\ar@{-->}[r]^g_{!}&Gx\ar[d]^{\psi_x}\\&x}\]
commute in $C$.

Now replace the above $F$ and $G$ with their restrictions to $D$ and
$E$, respectively.
Although
\cite[Section~IV.8, Theorem~1]{maclane} implies that
the restrictions
$F$ and $G$ give an adjoint pair of functors between the subcategories $D$ and $E$,
we want more: namely, we want
$F$ and $G$ to be
\emph{quasi-inverses}, so that
\[FG\cong \id_E\and GF\cong \id_D.\]
For this we impose a further
property on the universal arrows:
we assume that for $x\in \obj D$ and $y\in \obj E$ the
arrows $\eta_x$ and $\psi_y$ are actually ``universal on both
sides'', that is,
we assume that
not only is $(Fx,\eta_x)$ initial in $x\downarrow E$,
but also $(x,\eta_x)$ is final in $D\downarrow Fx$,
and similarly we assume that
not only is $(Gy,\psi_y)$ final in $D\downarrow y$,
but also $(y,\psi_y)$ is initial in $Gy\downarrow E$.
Summarizing, we assume:
\begin{enumerate}
\item for all $z,x\in \obj D$, $w\in \obj E$, $f\in C(z,Fx)$, and
$g\in C(x,w)$, there exist unique $h\in D(z,x)$ and $k\in E(Fx,w)$
such that
\[\xymatrix@C+30pt{
z \ar[dr]_f \ar@{-->}[r]^-h &x \ar[d]_{\eta_x} \ar[dr]^g
\\&Fx \ar@{-->}[r]_-k &w
}\]
commutes, and

\item for all $z\in \obj D$, $y,w\in \obj E$, $f\in C(z,y)$, and
$g\in C(Gy,w)$, there exist unique $h\in D(z,Gy)$ and $k\in E(y,w)$
such that
\[\xymatrix@C+30pt{
z\ar[dr]_f\ar@{-->}[r]^-h&Gy\ar[d]_{\psi_y}\ar[dr]^g
\\&y\ar@{-->}[r]_-k&w}\]
commutes.
\end{enumerate}

We include the elementary proof showing that the
above assumptions imply that $F$ and $G$ are quasi-inverses:

\begin{prop}
\label{subcategory}
With
the above notation and assumptions,
\[\xymatrix{D\ar@<2pt>[r]^F&E\ar@<2pt>[l]^G}\]
is a category equivalence.
\end{prop}

\begin{proof}
We aim to apply
\corref{equivalence2}, and for this
we need to show that $FG\cong \id_E$ and $F$
is full and faithful.
First, for each $x\in \obj E$, both
$(Gx,\psi_x)$ and $(FGx,\eta_{Gx})$ are initial in $Gx\dn E$,
so there is a unique isomorphism $\t_x$ in $E$ making the diagram
\[
\xymatrix{
&Gx \ar[dl]_{\psi_x} \ar[dr]^{\eta_{Gx}}
\\
x \ar[rr]_{\t_x}^\cong
&& FGx
}
\]
commute in $C$. We must show $\t$ is natural. Let $f:x\to y$ in $E$,
and consider the diagram
\[
\xymatrix{
x \ar[rr]^{\t_x}_\cong \ar[ddd]_f
&& FGx \ar[ddd]^{FGf}
\\
&Gx \ar[ul]^{\psi_x} \ar[d]_{Gf} \ar[ur]_{\eta_{Gx}}
\\
&Gy \ar[dl]_{\psi_y} \ar[dr]^{\eta_{Gy}}
\\
y \ar[rr]_{\t_y}^\cong
&& FGy.
}
\]
We need to show that the outer rectangle commutes, that is,
\[\t_y\circ f=FGf\circ \t_x.\]
It suffices to show that
\[\t_y\circ f\circ \psi_x=FGf\circ \t_x\circ \psi_x,\]
since $(x,\psi_x)$ is initial in $Gx\dn E$.
The left and right quadrilaterals and the top and bottom triangles all commute, so we have
\begin{align*}
\t_y\circ f\circ \psi_x
&=\t_y\circ \psi_y\circ Gf
\\&=\eta_{Gy}\circ Gf
\\&=FGf\circ \eta_{Gx}
\\&=FGf\circ \t_x\circ \psi_x.
\end{align*}
Thus $\t:FG\iso \id_E$.

Second, if $x,y\in \obj D$ and $g\in E(Fx,Fy)$, then since $(y,\eta_y)$
is final in $D\downarrow Fy$
there exists a unique $f\in D(x,y)$ such that
\[\xymatrix{
x\ar@{-->}[r]^-f\ar[d]_{\eta_x}
&y\ar[d]^{\eta_y}
\\
Fx\ar[r]_-g
&Fy
}\]
commutes. By definition of $F$, $f$ is the unique element of
$D(x,y)$ such that $g=Ff$.
Thus $F$ is full and faithful.
\end{proof}

Next we verify that the above equivalence carries over to
comma categories:
Let $a\in \obj D$,
and consider the subcategories $a\downarrow D$ and $a\downarrow E$
of the comma category $a\downarrow C$.

\begin{cor}
\label{comma}
Let $F:D\to E$ be an equivalence between subcategories of $C$ as in
\propref{subcategory}, and let $a\in \obj D$. Then the mapping
$\wilde F:a\downarrow D\to a\downarrow E$ defined by
\begin{align*}
\wilde F(x,f)&=(Fx,\eta_x\circ f)=(Fx,Ff\circ \eta_a)\\
\wilde Fh&=Fh\qquad
\text{for an arrow $h$ of $a\downarrow D$}
\end{align*}
is also an equivalence.
\end{cor}

\begin{proof}
It is routine to verify that $\wilde F$ is a functor.
We first show that $\wilde F$ is essentially surjective,
that is, that every object $(y,g)$ of $a\dn E$ is
isomorphic to one in the image of $\wilde F$.
Consider the diagram
\[\xymatrix{
a \ar@{-->}[r]^-h \ar[dr]_g
&Gy \ar[d]^{\psi_y} \ar[dr]^{\eta_{Gy}}
\\
&y \ar@{-->}[r]_-\t^-\cong
&FGy,
}\]
where $G:E\to D$ is the quasi-inverse functor from \corref{equivalence2}.
Since $(Gy,\psi_y)$ is
final in $D\downarrow y$, there exists $h\in D(a,Gy)$ such that
$g=\psi_y\circ h$. Since both $(y,\psi_y)$ and $(FGy,\eta_{Gy})$ are
initial in $Gy\downarrow E$, there is an isomorphism $\t:y\iso
FGy$ in $E$ such that $\t\circ \psi_y=\eta_{Gy}$. Then
\[\t:(y,g)\iso (FGy,\eta_{Gy}\circ h)=\wilde F(Gy,h)
\textin a\downarrow E.\]

Next we show that $\wilde F$ is full and faithful. Let
\[h:\wilde F(x,f)\to \wilde F(y,g)\textin a\downarrow E.\]
We must show that there exists a unique $k:(x,f)\to (y,g)$ in
$a\downarrow D$ such that $h=\wilde Fk$.
Consider the diagram
\[\xymatrix{
&a \ar[dl]_f \ar[dr]^g \ar[dd]^(.7){\eta_a}
\\
x \ar[dd]_{\eta_x} \ar@{-->}[rr]_(.3)k
&&y \ar[dd]^{\eta_y}
\\
&Fa \ar[dl]_{Ff} \ar[dr]^{Fg}
\\
Fx \ar[rr]_h
&&Fy
}\]
We want a unique $k$ for which the top triangle commutes and $h=\wilde
Fk=Fk$.
Since $F$ is full and faithful, there exists a unique $k\in D(x,y)$
such that $h=Fk$.
It suffices to show that the top triangle commutes, for then
$k$ will be the unique arrow from $(x,f)$ to $(y,g)$ in $a\downarrow D$
such that $h=\wilde Fk$.
Since $F$ is an equivalence, it further suffices to show that the bottom
triangle commutes.
Since $h$ is an arrow from $(Fx,\eta_x\circ f)=(Fx,Ff\circ \eta_a)$ to
$(Fy,\eta_y\circ g)=(Fy,Fg\circ \eta_a)$ in $a\downarrow E$, we have
\[h\circ Ff\circ \eta_a=\eta_y\circ g.\]
Since $Fg$ is the unique element of $E(Fa,Fy)$ such that $Fg\circ \eta_a=\eta_y\circ g$,
we must have
\[h\circ Ff=Fg,\]
and we are done.
\end{proof}

\section{Maximal and normal coactions}
\label{maximal normal}

Here we apply \propref{subcategory} and \corref{comma} to maximal
and normal coactions.

Maximalizations and normalizations have the
following universal properties (for maximalizations see
\cite[Lemma~6.2]{fischer:quantum}, and for normalizations see
\cite[Lemma~4.2]{fischer:quantum} and also
\cite[Lemma~2.1]{ekq} --- the requirement
in \cite{ekq}
that the homomorphisms are into the $C^*$-algebras themselves rather
than into multipliers is not used in the proof of \cite[Lemma~2.1]{ekq}):
\begin{enumerate}
\item
The subcategory $\CC^m(G)$ of $\CC(G)$ is coreflective,
and every maximalizing map $(B,\e)\to (A,\d)$ is a universal arrow from the inclusion functor $\CC^m(G)\hookrightarrow \CC(G)$ to $(A,\d)$.
In other words,
if $(B,\e,\p)$ is a maximalization of $(A,\d)$, $(C,\g)$ is a
maximal coaction, and $\tau:(C,\g)\to (A,\d)$ in $\CC(G)$, then
there is a unique arrow $\s$ making the diagram
\[\xymatrix{
(C,\g) \ar@{-->}[r]^-\s_-{!} \ar[dr]_\tau
&(B,\e) \ar[d]^\p
\\&(A,\d)
}\]
commute in $\CC(G)$.
\item
The subcategory $\CC^n(G)$ of $\CC(G)$ is reflective,
and every normalizing map $(A,\d)\to (B,\e)$ is a universal arrow from
$(A,\d)$ to
the inclusion functor $\CC^n(G)\hookrightarrow \CC(G)$.
In other words,
if $(B,\e,\p)$ is a normalization of $(A,\d)$, $(C,\g)$ is a
normal coaction, and $\tau:(A,\d)\to (C,\g)$ in $\CC(G)$, then
there is a unique arrow $\s$ such that
\[\xymatrix{
(A,\d) \ar[d]_\p \ar[dr]^\tau
\\(B,\e) \ar@{-->}[r]_-\s^-{!} &(C,\g)
}\]
commutes.
\end{enumerate}

The above universal properties have the following consequences:
as we observed in the abstract setting of~\secref{nonsense},
once a normalization $(A^n,\d^n,\eta_{(A,\d)})$ has been chosen for
every coaction $(A,\d)$,
there is a
unique functor from $\CC(G)$ to $\CC^n(G)$ given by
\[
(A,\d)\mapsto (A^n,\d^n)
\and
\p\mapsto \p^n,
\]
such that for every arrow $\p:(A,\d)\to (B,\e)$ in $\CC(G)$
the diagram
\[\xymatrix{
(A,\d) \ar[r]^-\p \ar[d]_{\eta_{(A,\d)}}
&(B,\e) \ar[d]^{\eta_{(B,\d)}}
\\
(A^n,\d^n) \ar[r]_-{\p^n}
&(B^n,\e^n)
}\]
commutes.
Similarly, once a maximalization $(A^m,\d^m,\psi_{(A,\d)})$ has been
chosen for
every coaction $(A,\d)$,
there is a
unique functor from $\CC(G)$ to $\CC^m(G)$ given by
\[
(A,\d)\mapsto (A^m,\d^m)
\and
\p\mapsto \p^m,
\]
such that for every arrow $\p:(A,\d)\to (B,\e)$ in $\CC(G)$
the diagram
\[\xymatrix{
(A^m,\d^m) \ar[r]^-{\p^m} \ar[d]_{\psi_{(A,\d)}}
&(B^m,\e^m) \ar[d]^{\psi_{(B,\d)}}
\\
(A,\d) \ar[r]_-\p
&(B,\e)
}\]
commutes.

The following proposition (further) justifies the terminology ``maximal'', and suggests that normal coactions could also be regarded as ``minimal'':

\begin{prop}
\label{characterize}
For any coaction $(A,\d)$,
\begin{enumerate}
\item $(A,\d)$ is maximal if and only if
$\p:(B,\e)\to (A,\d)$
is an isomorphism in $\CC(G)$
whenever
$\p\times G:B\times_\e G\to A\times_\d G$ is an isomorphism in $\CC$;

\item $(A,\d)$ is normal if and only if
$\p:(A,\d)\to (B,\e)$
is an isomorphism in $\CC(G)$
whenever
$\p\times G:A\times_\d G\to B\times_\e G$ is an isomorphism in $\CC$;
\end{enumerate}
\end{prop}

\begin{proof}
(i)
First assume $(A,\d)$ is maximal.
Let $\p:(B,\e)\to (A,\d)$ in $\CC(G)$, and assume that $\p\times G$ is an isomorphism.
Since the diagram
\[\xymatrix@C+60pt{
B\times_\e G\times_{\what\e} G
\ar[r]^-{\P_B}
\ar[d]_{\p\times G\times G}^\cong
&B\otimes \KK(L^2(G))
\ar[d]^{\p\otimes\id}
\\
A\times_\d G\times_{\what\d} G
\ar[r]_-{\P_A}^-\cong
&A\otimes \KK(L^2(G))
}\]
commutes,
$(\p\otimes\id)\circ \P_B$ is injective,
hence so is
$\P_B$.
Since $\P_B$ is also surjective, it is an isomorphism.
Therefore
\[\p\otimes\id=\P_A\circ (\p\times G\times G)\circ \P_B\inv\]
is also an isomorphism, hence so is $\p$.

Conversely, assume that
$\p:(B,\e)\to (A,\d)$
is an isomorphism
whenever
$\p\times G:B\times_\e G\to A\times_\d G$ is.
In particular, the maximalization map $\psi_{(A,\d)}:(A^m,\d^m)\to (A,\d)$ is an isomorphism since $\psi_{(A,\d)}\times G$ is. Therefore $(A,\d)$ is maximal since $(A^m,\d^m)$ is.

(ii)
First assume $(A,\d)$ is normal.
Let $\p:(A,\d)\to (B,\e)$ in $\CC(G)$, and assume that $\p\times G$ is an isomorphism.
Since the diagram
\[\xymatrix@C+60pt{
A \ar[r]^-{j_A}_-{\txt{injective}} \ar[d]_\p
&M(A\times_\d G) \ar[d]^{\p\times G}_\cong
\\
B \ar[r]_-{j_B}
&M(B\times_\e G)
}\]
commutes, $j_B\circ \p$ is injective, hence so is $\p$.
Since $\p$ is also surjective by the elementary \lemref{onto} below, it
is an isomorphism.

Conversely, assume that
$\p:(A,\d)\to (B,\e)$
is an isomorphism
whenever
$\p\times G:A\times_\d G\to B\times_\e G$ is.
In particular, the normalization map $\eta_{(A,\d)}:(A,\d)\to (A^n,\d^n)$ is an isomorphism since $\eta_{(A,\d)}\times G$ is. Therefore $(A,\d)$ is normal since $(A^n,\d^n)$ is.
\end{proof}

The above proof of (ii) used the following result,
which we state separately since it may be useful elsewhere:

\begin{lem}
\label{onto}
Let $\p:(A,\d)\to (B,\e)$ be an arrow in $\CC(G)$ such that
\[\p\times G:A\times_\d G\to B\times_\e G\]
is an isomorphism. Then $\p(A)=B$.
\end{lem}

\begin{proof}
Applying the maximalization functor gives a
commuting diagram
\begin{equation}
\label{commute}
\xymatrix@C+20pt{
(A^m,\d^m)
\ar[r]^-{\p^m}
\ar[d]_{\psi_{(A,\d)}}
&(B^m,\e^m)
\ar[d]^{\psi_{(B,\e)}}
\\
(A,\d)
\ar[r]^-\p
&(B,\e),
}
\end{equation}
and then the crossed-product functor and properties of maximalization give the commuting diagram
\[
\xymatrix@C+30pt{
(A^m\times_{\d^m} G,\what{\d^m},j_G)
\ar[r]^-{\p^m\times G}
\ar[d]_{\psi_{(A,\d)}\times G}^\cong
&(B^m\times_{\e^m} G,\what{\e^m},j_G)
\ar[d]_\cong^{\psi_{(B,\e)}\times G}
\\
(A\times_\d G,\what\d,j_G)
\ar[r]_-{\p\times G}^-\cong
&(B\times_\e G,\what\e,j_G).
}
\]
Thus $\p^m\times G$ is an isomorphism, hence so is $\p^m$,
by \corref{characterize}~(i).
In particular, $\p^m$ is surjective.
Since $\psi_{(B,\e)}$ is also surjective, so is $\p$, by commutativity of \diagref{commute}.
\end{proof}

We now show that the above universal properties also imply that maximal and normal coactions are equivalent:

\begin{thm}
\label{maximalnormalnocomma}
Normalization and maximalization give a category equivalence
\[\CC^m(G)\sim \CC^n(G).\]
\end{thm}

\begin{proof}
Observe that if $(A,\d)$ is a maximal coaction,
then $\eta_{(A,\d)}:(A,\d)\to (A^n,\d^n)$ is not only a normalization of $(A,\d)$, but also a maximalization of $(A^n,\d^n)$.
Similarly, if $(A,\d)$ is normal
then $\psi_{(A,\d)}:(A^m,\d^m)\to (A,\d)$ is not only a maximalization of $(A,\d)$, but also a normalization of $(A^m,\d^m)$.
Thus the normalizing and maximalizing maps are
``universal on both sides'', as in the hypotheses of
\propref{subcategory}. Therefore, that proposition tells us that
the restrictions of normalization to $\CC^m(G)$ and
maximalization to $\CC^n(G)$ give an equivalence $\CC^m(G)\sim \CC^n(G)$.
\end{proof}

As we have seen abstractly in \secref{nonsense}, the extension of all
the above to comma
categories is now automatic. We have in mind to use comma categories
of objects under the coaction
$(C^*(G),\d_G)$:

\begin{cor}
\label{maximalnormal}
Normalization gives an equivalence
\[(C^*(G),\d_G)\downarrow \CC^m(G)\sim (C^*(G),\d_G)\downarrow
\CC^n(G).\]
\end{cor}

\begin{proof}
Since the coaction $(C^*(G),\d_G)$ is maximal, we are precisely in
the situation of \propref{comma}, and the result follows immediately.
\end{proof}

Note that the normalization of an object $(A,\d,u)$ of
$(C^*(G),\d_G)\downarrow \CC^m(G)$ is
\[(A^n,\d^n,\eta_{(A,\d)}\circ u),\]
where $(A^n,\d^n,\eta_{(A,\d)})$ is the normalization of the coaction
$(A,\d)$.

\section{Actions and normal coactions}
\label{action normal}

Here we aim to show that the reduced-crossed-product functor gives
rise to a category equivalence. Thus, what the second author has
called ``Landstad duality'' \cite{Q92} is actually a category
equivalence.

\begin{thm}
\label{actionnormal}
The functor
\begin{align*}
(A,\a)&\mapsto (A\times_{\a,r} G,\what\a^r,i^r_G)
\\
\p&\mapsto \p\times_r G
\end{align*}
from $\AA(G)$ to the comma category $(C^*(G),\d_G)\downarrow
\CC^n(G)$
is an equivalence.
\end{thm}

Note that the definition of the functor is valid because
the coaction $\what\a^r$ is normal.

\begin{proof}
We need to verify that the reduced-crossed-product functor is full, faithful, and essentially surjective.
Essential surjectivity means that every object in $(C^*(G),\d_G)\dn \CC^n(G)$ is isomorphic to one of the form $(A\times_{\a,r} G,\what\a^r,i^r_G)$ for some action $(A,\a)$,
and is proved in
Landstad's characterization of reduced crossed products \cite[Theorem~3]{lan:dual}
(when translated from reduced to normal coactions, as in the discussion on page~\pageref{boilerplate})

To prove the functor is full and faithful,
suppose $(A,\a)$ and $(B,\b)$ are actions and
\[\s:(A\times_{\a,r} G,\what\a^r,i^{\a,r}_G)\to
(B\times_{\b,r} G,\what\b^n,i^{\b,r}_G)\]
is an arrow in $(C^*(G),\d_G)\downarrow \CC^n(G)$. We must show
that there exists a
unique arrow
\[\p:(A,\a)\to (B,\b)\]
in $\AA(G)$ such that $\s=\p\times_r G$.

Consider the diagram
\[\xymatrix{
A \ar[r]^-{i^r_A} \ar@{-->}[d]^?_\p
&M(A\times_{\a,r} G) \ar[d]^{\s}
\\
M(B) \ar[r]_-{i^r_B}
&M(B\times_{\b,r} G).
}\]

\subsection*{Step 1}
We first show that there is a unique homomorphism $\p$ making the
diagram
commute.
Since the homomorphism $i^r_B$ is faithful, it
suffices
to show that
\[\s\circ i^r_A(A)\subset i^r_B(M(B)).\]

Let $a\in A$. To show that $x:=\s\circ i^r_A(a)\in i^r_B(M(B))$,
by nondegeneracy of $i^r_B$  it suffices to show that
$x$ idealizes $i^r_B(B)$. Let $y\in i^r_B(B)$. We must show
that
\[xy,yx\in i^r(B).\]
Upon taking adjoints, it suffices to consider $xy$, and for this it
suffices to verify Landstad's conditions
\cite[Equations~(3.6)--(3.8)]{lan:dual}:
\newcommand{\landstad}[1]{(L\ref{#1})}
\begin{enumerate}
\renewcommand{\labelenumi}{(L\arabic{enumi})}
\item
\label{landstad one}
$\what\b^n(xy)=xy\otimes 1$,
\item
\label{landstad two}
$xyi^r_G(f),i^r_G(f)xy\in B\times_{\b,r} G$ for all
$f\in C_c(G)$, and
\item
\label{landstad three}
$s\mapsto \ad i^r_G(s)(xy)$ is norm continuous.
\end{enumerate}
Of course, $y$ itself satisfies these conditions.
Since $\s$ is $\what\a^r-\what\b^n$ equivariant, we have
\[\what\b^n\circ\s\circ i^r_A=\s\circ i^r_A\otimes 1,\]
which implies
\landstad{landstad one}
because $\what\b^n(y)=y\otimes 1$.

The first part of
\landstad{landstad two}
is immediate because $yi^r_G(f)\in
B\times_{\b,r}
G$. For the second part of
\landstad{landstad two},
let $y=i^r_B(b)$ with $b\in B$. Then
\begin{align*}
i^{\b,r}_G(f)\s\circ i^r_A(a)i^r_B(b)
&=\s\circ i^{\a,r}_G(f)\s\circ i^r_A(a)i^r_B(b)
\\&=\s\bigl(i^{\a,r}_G(f)i^r_A(a)\bigr)i^r_B(b)
\\&\approx
\s\left(\sum_{i=1}^ki^r_A(a_i)i^{\a,r}_G(f_i)\right)i^r_B(b)
\\&\qquad\righttext{for some}a_i\in A,f_i\in C_c(G)
\\&=\sum_{i=1}^k\s\circ i^r_A(a_i)\s\circ i^{\a,r}_G(f_i)i^r_B(b)
\\&=\sum_{i=1}^k\s\circ i^r_A(a_i)i^{\b,r}_G(f_i)i^r_B(b),
\end{align*}
which is in $B\times_{\b,r}G$ because each $i^{\b,r}_G(f_i)i^r_B(b)$
is.

For
\landstad{landstad three},
note that
\begin{align*}
\ad i^{\b,r}_G(s)(x)
&=\ad \s\circ i^{\a,r}_G(s)\bigl(\s\circ i^r_A(a)\bigr)
\\&=\s\bigl(\ad i^{\a,r}_G(s)\circ i^r_A(a)\bigr)
\\&=\s\bigl(i^r_A\circ\a_s(a)\bigr),
\end{align*}
which is norm continuous, and
\landstad{landstad three}
follows because $s\mapsto \ad
i^{\b,r}_G(s)(y)$ is also norm continuous.

Thus $\s:i^r_A(A)\to M(i^r_B(B))=i^r_B(M(B))$, and by injectivity of
$i^r_B$ there is a
unique homomorphism $\p:A\to M(B)$ such that $i^r_B\circ \p=\s\circ
i^r_A$.

\subsection*{Step 2}

We show that $\p:A\to M(B)$ is $\a-\b$ equivariant: again by
injectivity of $i^r_B$ the following computation suffices:
\begin{align*}
i^r_B\circ\p\circ\a_s
&=\s\circ i^r_A\circ\a_s
\\&=\s\circ\ad i^{\a,r}_G(s)\circ i^r_A
\\&=\ad \s(i^{\a,r}_G(s))\circ\s\circ i^r_A
\\&=\ad i^{\b,r}_G(s)\circ i^r_B\circ\p
\\&=i^r_B\circ\b_s\circ\p.
\end{align*}

\subsection*{Step 3}
We show that $\p:A\to M(B)$ is nondegenerate.
Suppose not. Then we can find a nonzero functional
$\psi\in B^*$ which annihilates $\p(A)B$.
Our strategy is to use crossed-product duality: the
double-crossed-product homomorphism
\[\p\times_r G\times G:A\times_{\a,r} G\times_{\what\a^r} G
\to
M(B\times_{\b,r} G\times_{\what\b^n} G)\]
is nondegenerate because $\p\times_r G$ is, and by
the Imai-Takai duality theorem (see \cite[Theorem~A.67]{enchilada})
we have an isomorphism
\[\Phi=(\id\otimes M)\circ \wilde \b\times
(1\otimes\l)\times (1\otimes M):
B\times_{\b,r} G\times_{\what\b^n} G\iso B\otimes \KK(L^2(G)),\]
where
\[\wilde\b:B\to \wilde M(B\otimes C_0(G))=C_b(G,B)\]
is defined by
\[\wilde \b(b)(s)=\b_s\inv(b)\for b\in B,s\in G.\]
Let $\xi,\eta\in L^2(G)$ be nonzero. Then the functional
\[\psi\otimes \omega_{\xi,\eta}\in (B\otimes \KK(L^2(G)))^*\]
is nonzero, where $\omega_{\xi,\eta}$ is the vector functional
\[\<T,\omega_{\xi,\eta}\>=(T\xi,\eta)\for T\in \KK(L^2(G)).\]
The composition
\[\pi:=\Phi\circ (\p\times_r G\times G)\circ j_{A\times G}\circ i^r_A:
A\to M(B\otimes \KK(L^2(G)))\]
is nondegenerate, and so is the representation $M:C_0(G)\to B(L^2(G))$, so we can find $a\in A$, $b\in B$, $f\in C_0(G)$, and $T\in \KK(L^2(G))$ such that
\[\<\pi(a)(b\otimes M_fT),\psi\otimes \omega_{\xi,\eta}\>\ne 0.\]
We have
\[(\p\times_r G\times G)\circ j_{A\times G}\circ i^r_A
=j_{B\times G}\circ (\p\times_r G)\circ i^r_A
=j_{B\times G}\circ i^r_B\circ \p,\]
so
\[\pi=(\id\otimes M)\circ \wilde\b\circ \p.\]
Also,
\[\wilde\b\circ \p(a)(b\otimes f)\in B\otimes C_0(G)=C_0(G,B)\]
is the function defined by
\begin{align*}
\bigl(\wilde\b\circ \p(a)(b\otimes f\bigr)(s)
&=\wilde\b\circ \p(a)(s)(b\otimes f)(s)
\\&=\b_s\inv(\p(a))bf(s)
\\&=\p(\a_s\inv(a))bf(s).
\end{align*}
Claim: if $g\in B\otimes C_0(G)$ and $\kappa\in L^2(G)$ then
\[\<(\id\otimes M)(g),\psi\otimes \omega_{\kappa,\eta}\>
=\int \<g(s),\psi\>\kappa(s)\bar{\eta(s)}\,ds.\]
This is easily checked on elementary tensors $b=c\otimes h$ with $c\in B$ and $h\in C_0(G)$:
\begin{align*}
\<(\id\otimes M)(c\otimes h),\psi\otimes \omega_{\kappa,\eta}\>
&=\<c\otimes M_h,\psi\otimes \omega_{\kappa,\eta}\>
\\&=\<c,\psi\>\<M_h,\omega_{\kappa,\eta}\>
\\&=\<c,\psi\>(M_h\kappa,\eta)
\\&=\<c,\psi\>\int h(s)\kappa(s)\bar{\eta(s)}\,ds
\\&=\int \<c,\psi\>h(s)\kappa(s)\bar{\eta(s)}\,ds
\\&=\int \<ch(s),\psi\>\kappa(s)\bar{\eta(s)}\,ds
\\&=\int \<(c\otimes h)(s),\psi\>\kappa(s)\bar{\eta(s)}\,ds,
\end{align*}
and the claim follows by linearity, density, and continuity. Thus
\begin{align*}
0
&\ne \<\pi(a)(b\otimes M_fT),\psi\otimes \omega_{\xi,\eta}\>
\\&=\<(\id\otimes M)\circ \wilde\b\circ \p(a)(\id\otimes M)(b\otimes f),
\psi\otimes \omega_{T\xi,\eta}\>
\\&=\<(\id\otimes M)\bigl(\wilde\b\circ \p(a)(b\otimes f)\bigr),
\psi\otimes \omega_{T\xi,\eta}\>
\\&=\int \<\bigl(\wilde\b\circ \p(a)(b\otimes f)\bigr)(s),\psi\>
(T\xi)(s)\bar{\eta(s)}\,ds
\\&=\int \<\p(\a_s\inv(a))bf(s),\psi\>
(T\xi)(s)\bar{\eta(s)}\,ds
\\&=\int \<\p(\a_s\inv(a))b,\psi\>f(s)
(T\xi)(s)\bar{\eta(s)}\,ds
=0,
\end{align*}
because
\[\<\p(\a_s\inv(a)b,\psi\>=0\all s\in G,\]
and we have a contradiction.

\subsection*{Step 4}
Finally, we verify that $\s=\p\times_r G$: for $a\in A,f\in C_0(G)$ we have
\begin{align*}
\p\times_r G\bigl(i^r_A(a)i^r_G(f)\bigr)
&=i^r_B(\p(a))i^r_G(f)
=\s\bigl(i^r_A(a)i^r_G(f)\bigr).
\qedhere
\end{align*}
\end{proof}

\begin{rem}
\label{preserve}
Like any category equivalence, the reduced-crossed product functor preserves some properties of morphisms. For example, if $\p:(A,\a)\to (B,\b)$ in $\AA(G)$, then $\p:A\to B$ is an isomorphism in $\CC$ if and only if $\p\times_r G:A\times_{\a,r} G\to B\times_{\b,r} G$ is --- just note that $\p\times_r G$ is an isomorphism in $\CC$ if and only if it is an isomorphism in the comma category $(C^*(G),\d_G)\dn \CC^n(G)$.
\end{rem}

\section{Actions and maximal coactions}
\label{action maximal}

Here we combine Theorems~\ref{maximalnormal} and \ref{actionnormal} to
show that the full-crossed-product functor gives a category
equivalence.

\begin{thm}
\label{full landstad}
The functor
\begin{align*}
(A,\a)&\mapsto (A\times_\a G,\what\a,i_G)
\\
\p&\mapsto \p\times G
\end{align*}
from $\AA(G)$ to the comma category $(C^*(G),\d_G)\downarrow \CC^m(G)$
is an equivalence.
\end{thm}

\begin{proof}
\thmref{actionnormal} says that the reduced-crossed-product functor
gives an equivalence
\[\AA(G)\sim (C^*(G),\d_G)\downarrow \CC^n(G),\]
and \corref{maximalnormal} says that normalization gives an
equivalence
\[(C^*(G),\d_G)\downarrow\CC^m(G)\sim
(C^*(G),\d_G)\downarrow\CC^n(G).\]
Thus it suffices to show that
the composition
\[(A,\a)\mapsto (A\times_\a G,\what\a,i_G)\mapsto
((A\times_\a G)^n,(\what\a)^n,\eta_{(A\times_\a G,\what\a)}\circ
i_G),\]
that is, full crossed product followed by normalization,
is naturally isomorphic to the reduced-crossed-product functor
\[(A,\a)\mapsto (A\times_{\a,r} G,\what\a^r,i^r_G).\]
Thus, for each action $(A,\a)$ we need an isomorphism
\[(A\times_{\a,r} G,\what\a^r,i^r_G)\cong
((A\times_\a G)^n,(\what\a)^n,\eta_{(A\times_\a G,\what\a)}\circ
i_G)\]
which is natural in $(A,\a)$.
Since both
\[(A\times_{\a,r} G,\what\a^r,\L_{(A,\a)})
\and
((A\times_\a G)^n,(\what\a)^n,\eta_{(A\times_\a G,\what\a)})\]
are normalizations of $(A\times_\a G,\what\a)$,
there is a unique isomorphism $\t_{(A,\a)}$ in $\CC^n(G)$ such that
\[\xymatrix{
(A\times_{\a,r} G,\what\a^r)
\ar@{-->}[rr]^-{\t_{(A,\a)}}_-\cong
&&((A\times_\a G)^n,(\what\a)^n)
\\
&(A\times_\a G,\what\a)
\ar[ul]^{\L_{(A,\a)}}
\ar[ur]_{\eta_{(A\times_\a G,\what\a)}}
}\]
commutes in $\CC(G)$.
We need to know that the above commuting triangle
passes over to one in the comma category $(C^*(G),\d_G)\downarrow
\CC(G)$:
\[\xymatrix@C=0in{
(A\times_{\a,r} G,\what\a^r,i^r_G)
\ar[rr]^-{\t_{(A,\a)}}
&&((A\times_\a G)^n,(\what\a)^n,
\eta_{(A\times_\a G,\what\a)}\circ i_G)
\\
&(A\times_\a G,\what\a,i_G)
\ar[ul]^{\L_{(A,\a)}}
\ar[ur]_{\eta_{(A\times_\a G,\what\a)}}
}\]
and for this it suffices to show that
\[\t_{(A,\a)}\circ i^r_G=\eta_{(A\times_\a G,\what\a)}\circ i_G.\]
Consider the diagram
\[\xymatrix{
A\times_{\a,r} G \ar[rr]^-{\t_{(A,\a)}}
&&A\times_\a G
\\
&C^*(G)
\ar[ul]_{i^r_G} \ar[ur]^{i_G} \ar[d]^(.4){i_G}
\\
&A\times_\a G
\ar[uul]^{\L_{(A,\a)}}
\ar[uur]_{\eta_{(A\times_\a G,\what\a)}}.
}\]
The outer triangle commutes by the above,
the left triangle commutes by definition of $i^r_G$,
and the right triangle commutes because
$\eta_{(A\times_\a G,\what\a)}$ is an arrow in $(C^*(G),\d_G)\downarrow
\CC(G)$.
Therefore the top triangle commutes, as desired.

It remains to verify that the isomorphisms $\t$ are natural:
Let
\[\p:(A,\a)\to (B,\b)\textin\AA(G),\]
and consider the diagram
\[\xymatrix@C=0in{
(A\times_{\a,r} G,\what\a^r,i^r_G)
\ar[rr]^-{\t_{(A,\a)}}
\ar[ddd]_{\p\times_r G}
&&((A\times_\a G)^n,(\what\a)^n,\eta_{(A\times_\a G,\what\a)}\circ
i_G)
\ar[ddd]^{(\p\times G)^n}
\\
&(A\times_\a G,\what\a,i_G)
\ar[ul]^{\L_{(A,\a)}}
\ar[ur]_{\eta_{(A\times_\a G,\what\a)}}
\ar[d]^{\p\times G}
\\
&(B\times_\b G,\what\b,i_G)
\ar[dl]_{\L_{(B,\b)}}
\ar[dr]^{\eta_{(B\times_\b G,\what\b)}}
\\
(B\times_{\b,r} G,\what\b^n,i^r_G)
\ar[rr]_-{\t_{(B,\b)}}
&&((B\times_\b G)^n,(\what\b)^n,\eta_{(B\times_\b G,\what\b)}\circ
i_G).
}\]
The top and bottom triangles commute by definition of $\t$, the left
quadrilateral commutes by naturality of $\L$, and the right
quadrilateral commutes by definition of $(\p\times G)^n$. Since
$\L_{(A,\a)}$ is surjective, the outer rectangle commutes.
\end{proof}

\begin{rem}
It is instructive to compare \thmref{full landstad} with
\cite[Theorem~3.2]{KQ07}, which states that a $C^*$-algebra
$A$ is isomorphic to a full crossed product $B\times_\a G$ if and only
if there are a maximal coaction $\d$ of $G$ on $A$ and a strictly
continuous unitary homomorphism $u:G\to M(A)$ such that
$\d(u_s)=u_s\otimes s$. Thus, in the terminology of the present paper,
\cite[Theorem~3.2]{KQ07} shows that for every object
$(A,\d,u)$ in $(C^*(G),\d_G)\downarrow \CC^m(G)$ the associated object
$(A,\d)$ of $\CC^m(G)$ is isomorphic to one of the form $(B\times_\a
G,\what\a)$. \cite[Theorem~3.2]{KQ07} also states that the
action $(B,\a)$ is uniquely determined up to isomorphism subject to
the condition that, if $\t:(A,\d)\cong (B\times_\a G,\what\a)$ is the
isomorphism, then
\[
\L\circ\t\circ u=i^r_G.
\]
\thmref{full landstad} above is stronger, not only because it contains categorical information, but more importantly because it lifts the uniqueness clause to
\[
\t\circ u=i_G.
\]
The
heart of the matter is Fischer's universal
property \cite[Lemma~6.2]{fischer:quantum}, which implies the following:

\begin{cor}
If $(B,G,\a)$ is an action, then
$i_G:G\to
M(B\times_\a G)$ is the only strictly continuous unitary
homomorphism such that
\[\what\a(i_G(s))=i_G(s)\otimes s
\and
\L\circ i_G=i^r_G.\]
\end{cor}
The special case where $B=\C$, namely:
the canonical embedding $u:G\to M(C^*(G))$ is the only strictly continuous unitary homomorphism such that
\[\d_G(u_s)=u_s\otimes s
\and
\l(u_s)=\l_s,\]
is elementary, and does not require the results of this paper.
\end{rem}

\begin{rem}
\remref{preserve} has an obvious and true analogue for the full-crossed product functor; for example, if $\p:(A,\a)\to (B,\b)$ in $\AA(G)$ then $\p:A\to B$ is an isomorphism in $\CC$ if and only if $\p\times G:A\times_\a G\to B\times_\b G$ is.
\end{rem}


\providecommand{\bysame}{\leavevmode\hbox to3em{\hrulefill}\thinspace}
\providecommand{\MR}{\relax\ifhmode\unskip\space\fi MR }
\providecommand{\MRhref}[2]{%
  \href{http://www.ams.org/mathscinet-getitem?mr=#1}{#2}
}
\providecommand{\href}[2]{#2}

\end{document}